\documentclass[a4paper,11pt]{amsart}
\usepackage{fullpage,amssymb,latexsym,amsmath,amsthm}

\newtheorem{theorem}{Theorem}
\newtheorem{prop}{Proposition}
\newtheorem{corollary}{Corollary}
\theoremstyle{definition} \newtheorem{definition}{Definition}
\theoremstyle{remark}
\newtheorem{remark}{Remark}

\renewcommand{\SS}{\mathbb{S}}
\renewcommand{\phi}{\varphi}
\newcommand{\HH}{\mathbb{H}}
\newcommand{\XX}{\mathbb{X}}
\newcommand{\Lf}{\mathcal{L}}
\newcommand{\Rf}{\mathcal{R}}

\newcommand{\be}{\begin{equation}}
\newcommand{\ee}{\end{equation}}
\newcommand{\la}{\langle}
\newcommand{\ra}{\rangle}
\newcommand{\ba}{\begin{aligned}}
\newcommand{\ea}{\end{aligned}}

\newcommand{\Bf}{\mathcal{B}}
\newcommand{\Ff}{\mathcal{F}}

\def\Re{\mathbb{R}}

\def\1{1\!\!\hbox{{\rm I}}}
\begin{document}

 \title{Lift zonoid and barycentric representation on a Banach space with a cylinder measure}
 \author{Alexei M. Kulik}
\address{Institute of Mathematics Nat. Acad. Sci. of Ukraine, Tereshchenkivska str. 3,  01601, Kyiv, Ukraine}
\email{kulik@imath.kiev.ua}

\author{Taras D. Tymoshkevych}
 \address{National Taras Shevchenko University of Kyiv, Academician Glushkov pr. 4-e, 03127, Kyiv, Ukraine.}
 \email{tymoschkevych@gmail.com}




 \begin{abstract}
We show that the lift zonoid concept for a probability measure on $\Re^d$,
introduced  in (Koshevoy and Mosler, 1997), leads naturally to  a
one-to one representation of any interior point of the convex hull of
the support of a continuous measure as the barycenter w.r.t. to this
measure of either of a half-space, or the whole space. We prove the
infinite-dimensional generalization of this representation, which is
based on the extension of the lift-zonoid concept for a cylindrical
probability measure.
 \end{abstract}
 \maketitle 
%
\section{INTRODUCTION}\label{s:1}
For a  probability measure $\mu$ with a finite first moment, defined
on the Borel $\sigma$-algebra in $\Re^d$, its \emph{zonoid} $Z(\mu)$
is defined as the set of all the points in $\Re^d$ of the form
\be\label{int_mu} \int_{\Re^d}g(x)x\, \mu(dx) \ee with arbitrary
measurable $g:\Re^d\to [0,1]$, see \cite{KM_Be}, Definition 2.1. The
\emph{lift zonoid} $\hat Z(\mu)$ is defined as the zonoid of the
measure $\delta_1\times\mu$ in $\Re^{d+1}$. Equivalently, the zonoid
$Z(\mu)$ and the lift zonoid $\hat Z(\mu)$ are the sets of the
points of the form \be\label{expect} Eg(X)X\in \Re^d\quad
\hbox{and}\quad \Big(Eg(X), Eg(X)X\Big)\in \Re^{d+1} \ee
respectively, where $X$ is a random vector with the distribution
$\mu$.  The zonoid $Z(\mu)$ (resp.,  the lift zonoid  $\hat Z(\mu)$)
is a  convex compact set in $\Re^d$ (resp., $\Re^{d+1}$), symmetric
w.r.t. the point $(1/2)EX$ (resp., the point $((1/2),(1/2)EX)$).

Zonoid and lift zonoid concepts appear to be  useful in various
aspects. The lift zonoid determines the underlying measure uniquely;
see \cite{KM_Be}, \cite{KM_AnnSta}. This  important observation
leads to natural applications of the lift zonoid concept in
multivariate statistics for measuring the variability of  laws of
random vectors, and for ordering these laws, see \cite{KM_Be}.
Although the zonoid of a measure does not determine the measure
uniquely,  the concept of \emph{zonoid equivalence} appears to be
both naturally motivated by financial applications, and useful for
proving extensions of the ergodic theorem for zonoid stationary and
zonoid swap-invariant random sequences, see \cite{MS11},
\cite{MSS12}. Lift zonoids also lead to  natural  definitions of
associated $\alpha$-trimming and data depth, see \cite{KM_AnnSta}.

In this paper, we are concentrated  on the application of the lift
zonoid concept  to a barycentric representation of the points of the
initial space w.r.t. the given measure $\mu$. The possibility of
such an application was mentioned in \cite{KM_AnnSta}; however,
respective Theorem 3.1 in \cite{KM_AnnSta} is not the best possible
one. Below we show that, for a \emph{continuous} measure $\mu$ (see
condition (\ref{lines}) below), every interior  point of the closed
convex support  $\mathbf{S}(\mu)$ of the measure $\mu$ has a
\emph{unique} representation as a barycenter either of a half-space,
or of the whole space. This  barycentric representation has the same
spirit with the classical Krein-Milman theorem, see e.g.
\cite{Rud73}, Chapter 3, but unlike the Krein-Milman theorem  has an
important uniqueness feature. In some particularly
interesting cases, this makes it possible to introduce the
\emph{barycentric coordinates} on  the initial space.

The lift zonoid concept can be extended naturally to the
infinite-dimensional setting, and such an extension  appears to be
particularly  useful. For instance, it was shown in \cite{Borell09}
that the (analogue of) the lift zonoid induced by the Brownian
motion can be used in the sensitivity analysis of a
certain iterative portfolio selection problem. In this paper, we
show that infinite-dimensional lift zonoids can be used efficiently
to obtain a barycentric representation on a Banach space. It is
important observation that, in order to obtain the  barycentric
representation  in its most closed and simple form, one needs to
consider a \emph{cylinder probability measure} $\mu$ on the initial
space, rather than a usual probability measure. This leads naturally
for the extension of  the lift zonoid concept for  cylinder
probability measures.

The structure of the paper is the following. In Section 2, we
introduce the notions of the zonoid, the lift zonoid, and the zonoid
trimmed regions, induced by a cylinder probability measure on a
Banach space. We give the basic properties of these objects, and
formulate the main result of this paper about the barycentric
representation w.r.t. to the given cylinder measure $\mu$. The
proofs of these results are given in Section 3. In Section 4, we
discuss the way how does the barycentric representation lead to three
closely related types of barycentric coordinates on the initial
space. The barycentric representation from Section 2 relies on the
assumption that \emph{the zonoid span} $\Lf(\mu)$ coincides with the
initial space; see the definition and the discussion in Remark
\ref{r3}. For a $\sigma$-additive measure $\mu$ such an assumption
typically fails; in Section 5 we formulate the version of the
barycentric representation theorem for a space with a
$\sigma$-additive measure $\mu$, which takes into account the
natural topology of the zonoid span $\Lf(\mu)$. In Section 6, we
consider an illustrative example with $\mu$ being either a centered
cylinder Gaussian measure with identity covariance operator, or a
centered $\sigma$-additive Gaussian measure.

\section{Main results}\label{s:2}

\subsection{Infinite-dimensional zonoid, lift zonoid, and zonoid trimmed regions, induced by a cylinder probability measure}

Let $\XX$ be a reflexive separable Banach space. Denote by
$\mathcal{C}$ its \emph{cylinder algebra}; that is, the family of
subsets in $\XX$  of the form
$$
\{x\in \XX: (\la x, x^*_1\ra, \dots \la x, x^*_m\ra)\in B\}, \quad
m\geq 1, \quad x_1^*, \dots, x_m^*\in \XX^*, \quad B\in \Bf(\Re^m),
$$
here and below   $\XX^*$ denotes the dual space to $\XX$, and
$\la\cdot, \cdot\ra$ denotes the duality between $\XX$ and $\XX^*$.
Let $\mu$ be a\emph{ cylinder probability measure} on $\mathcal{C}$
\emph{of type 1} (see \cite{VahTarCho}, Chapter IV, \S5), that is,
an additive function $\mu:\mathcal{C}\to [0,1]$ such that
\begin{itemize}\item[(a)] for every $m\geq1$, $x_1^*, \dots, x_m^*\in \XX^*$ the restriction of $\mu$ to the $\sigma$-algebra
$$
\mathcal{C}_{x_1^*, \dots, x_m^*}:=\Big\{\{x:(\la x, x^*_1\ra, \dots
\la x, x^*_m\ra)\in B\}, B\in \Bf(\Re^m)\Big\}=\sigma\Big(\la \cdot,
x^*_1\ra, \dots, \la \cdot, x^*_m\ra\Big)
$$
is a probability measure;
\item[(b)] there exists $C_\mu>0$ such that
\be\label{type1} \int_\XX|\la x, x^*\ra|\mu(dx)\leq
C_\mu\|x^*\|_{\XX^*}, \quad x^*\in \XX^*. \ee
\end{itemize}

The cylinder measure $\mu$ can be naturally interpreted as the
``law'' of a generalized random element $X$ in $\XX$ with finite
1-st order, which, by definition, is the linear bounded map
$$
X:\XX^*\to L_1(\Omega, \Ff, P),
$$
where $(\Omega, \Ff, P)$ is some probability space; see
\cite{VahTarCho}, Chapter IV for details.  The mean value $\int_\XX
x\mu(dx)$  or, equivalently, the expectation $E[X]$ is well defined
in the following sense. It can be seen easily that the map
$$
T:\XX^*\ni x^*\mapsto \int_\XX\la x, x^*\ra\mu(dx)=E\la X, x^*\ra\in
\Re
$$
is linear and continuous. Because $\XX$ is assumed to be reflexive,
there exists $x_\mu\in \XX$ such that $T(x^*)=\la x_\mu, x^*\ra$. By
the definition,
$$
\int_\XX x\mu(dx)=E[X]=x_\mu.
$$
Note that, by the same argument, for any bounded random variable
$\eta$  the expectation $E[\eta X]\in \XX$ is well defined.

\begin{definition} The zonoid $Z(\mu)$ of the cylinder measure $\mu$ is the set of all points in $\XX$ of the form $
E[\eta X]$ where $X$ is a generalized random element in $\XX$ with
the distribution $\mu$, and $\eta:\Omega\to [0,1]$ is arbitrary
$\sigma(X)$-measurable random variable.

The lift zonoid  $\hat Z(\mu)$ is  the set of all points in
$\Re\times \XX$ of the form $(E\eta, E[\eta X])$ with is arbitrary
$\sigma(X)$-measurable random variable $\eta:\Omega\to [0,1]$.

For $\alpha\in (0,1]$, the \emph{zonoid $\alpha$-trimmed region}
$D_\alpha(\mu)$  is the set of all the points in $\XX$ of the form
$$
{1\over \alpha}E[\eta X], \quad g\in F_\alpha,
$$
where $F_\alpha$ denotes the family of  $\sigma(X)$-measurable
random variables $\eta:\Omega\to [0,1]$ such that $E\eta=\alpha.$
\end{definition}

\begin{remark} If $X$ is a random vector in $\Re^d$, then any $\sigma(X)$-measurable random variable $\eta$ have the form $\eta=g(X)$ with a Borel measurable $g:\Re^d\to \Re$. Therefore, the above definition is a straightforward generalization of the finite-dimensional definitions of zonoid and  lift zonoid, mentioned in Introduction. \end{remark}

In the following theorem, the basic properties of   zonoid, lift
zonoid, and zonoid $\alpha$-trimmed regions are collected. Denote by
$\mathcal{H}^\XX$ the class of all sets of the form $\{x \in \XX:
\la x, x^*\ra\geq a\}$ with   $x^*\in \XX$,  $a\in \Re$; clearly,
every set $H\in \mathcal{H}^\XX$ is a cylinder set.  Denote also
$$
\mathbf{S}(\mu)=\bigcap_{H\in \mathcal{H}^\XX: \mu(H)=1} H;
$$
note that, in the case of $\XX=\Re^d$ and a (usual) probability
measure $\mu$, the set $\mathbf{S}(\mu)$  coincides with the closed
convex hull of the topological  support of $\mu$.

\begin{theorem}\label{t1} Let $\mu$ be a  cylinder probability measure of type 1, and $X$ be respective generalized random element.

Then the following properties hold true.
\begin{itemize}

\item[(a)] $Z(\mu)$ and $\hat Z(\mu)$ are bounded closed convex sets in $\XX$ and $\Re\times \XX$, respectively.

\item[(b)] $Z(\mu)$ and $\hat Z(\mu)$ are symmetric w.r.t. the points $(1/2)E[X]$ and $((1/2),(1/2)E[X])$, respectively.

\item[(c)] $\hat Z(\mu)$ determines $\mu$ uniquely.

\item[(d)] Every $D_\alpha(\mu), \alpha\in (0,1]$ is a bounded closed convex compact set in $\Re^d$.

\item[(e)] $D_\alpha(\mu) \supset D_\beta(\mu)$ for every $\alpha < \beta$.

\item[(f)] If the measure $\mu$ is centered (that is, $E[X]=0$), then
$$
\alpha D_\alpha(\mu)=-(1-\alpha)D_{1-\alpha}(\mu),\quad
\alpha\in(0,1].
$$

\item[(g)] For any sequence $\alpha_n\to \alpha\in (0,1)$,
$$
D_{\alpha_n}(\mu)\to D_\alpha(\mu)
$$
in the Hausdorff distance. At the point $\alpha=1$ the family
$D_\alpha(\mu), \alpha\in (0,1]$ is continuous in the following
weaker sense:
$$
D_1(\mu)=\bigcap_{\alpha\in (0,1)}D_\alpha(\mu).
$$

\item[(h)] The closure of $D_0(\mu):=\bigcup_{\alpha\in (0,1]}D_\alpha(\mu)$ equals
 $\mathbf{S}(\mu)$.
\end{itemize}
\end{theorem}

Theorem \ref{t1} is proved in Section \ref{s3} below.

It was mentioned in \cite{Borell09} briefly that, for a (usual)
probability measure $\mu$ its zonoid is necessarily a compact. More
precisely, this property can be formulated as follows.

\begin{prop}\label{p1} Let $\mu$ a probability measure such that the respective family of random variables
\be\label{fam} \{\la X, x^*\ra, \|x^*\|_{\XX^*}\leq 1\} \ee is
uniformly integrable.

Then  $Z(\mu)$,  $\hat Z(\mu)$, and  $D_\alpha(\mu), \alpha\in
(0,1]$ are compact sets in $\XX$.
\end{prop}

Proposition \ref{p1} is proved in Section \ref{s3}.

\begin{remark} It is sufficient for the uniform integrability of the family (\ref{fam}) that either $\mu$ has a strong moment of order $1$; that is
\be\label{so} E\|X\|<\infty, \ee or $\mu$ has a weak moment of some
order $p>1$; that is, \be\label{wo} E|\la X, x^*\ra|^p<\infty, \quad
x^*\in \XX^*. \ee Under the assumption that $\mu$ has a weak moment
of order $p>1$, the compactness property, equivalent to that of
Proposition \ref{p1}, was proved in \cite{Dorogovtsev2000}.
\end{remark}

 In the finite-dimensional setting, an alternative definition of the zonoid is the expectation of the random segment $[0,X]$, see \cite{KM_Be}.  To complete the list of main properties of the zonoid of a cylinder probability measure, we give the extension of this definition for the infinite-dimensional setting.

Denote by $\SS_{\XX^*}$  the unit sphere in $\XX^*$, and by $h(C,
\cdot)$ the \emph{support function} of a convex set $C\subset \XX$:
$$
h(C, x^*)=\sup\{\la x,x^*\ra, x\in C\},\quad x^*\in \SS_{\XX^*}.
$$
Consider  the ``random segment'' $[0,X]$, which is generalized in
the sense that $X$ is a generalized random element in $\XX$. Note
that for every $x^*\in \SS_{\XX^*}$ respective value of the support
function
$$
h([0,X], x^*)=\la X, x^*\ra_+:=\Big(\la X, x^*\ra\vee 0\Big)
$$
is measurable w.r.t. $\mathcal{C}_{x^*} $, dominated by $|\la X,
x^*\ra|$, and therefore is integrable.
 Define \emph{the expectation} $E[0,X]$ of $[0,X]$ as the convex closed set in $\XX$ such that
\be\label{exp}
 h(E[0,X], x^*)=E h([0,X], x^*), \quad x^*\in \SS_{\XX^*}.
 \ee
\begin{prop}\label{p11} \begin{enumerate}\item The set $E([0,X])$ is well defined by (\ref{exp}), and $Z(\mu)=E([0,X])$.
\item $\hat Z(\mu)=E[0,(1,X)]$.
\end{enumerate}
\end{prop}

The proof of Proposition \ref{p11} is close to that of Proposition
2.2 in \cite{KM_Be}, with a minor difference caused by the fact
$Z(\mu)$ is not necessarily a compact now; see  Section \ref{s3}
below.

\subsection{Barycentric representation}\label{s22}

In this section, introduce the \emph{barycentric representation} in
$\XX$ w.r.t. to the given cylinder measure $\mu$. We assume that
$\mu$ is
 \emph{continuous} in  the sense that
\be\label{lines} \mu(\{ x : \la x, x^*\ra = a \})=0, \quad x^*\in
\XX\setminus\{0\}, \quad a\in \Re. \ee For a cylinder set $A$ with
$\mu(A)>0$, its \emph{barycenter} w.r.t. $\mu$ is defined by the
identity
$$B_A(\mu)=\frac{\int_{A} x \mu (dx)}{\mu(A)}$$
(another term for the barycenter $B_A(\mu)$ is \emph{$A$-centroid of
$\mu$}, see \cite{KM_AnnSta}).

\begin{theorem}\label{t2} Let $\mu$ be a  cylinder probability measure of type 1. Assume that

\begin{itemize}\item[(i)] $\mu$ satisfy (\ref{lines});
\item[(ii)] the interior of  $D_0(\mu)$ is non-empty.
\end{itemize}

Then the following holds true.

\begin{enumerate}\item The set $\mathbf{B}(\mu)=\{B_A(\mu), \mu(A)>0\}$ of all barycenters w.r.t. $\mu$ equals the interior of $\mathbf{S}(\mu)$.
\item Any point $x\in \mathbf{B}(\mu)$ has a representation in the form
\be\label{bary}x=B_H(\mu), \quad H\in \mathcal{H}^\XX, \quad
\mu(H)>0.\ee This representation is unique in the following sense:
if  $B_{H}(\mu)=B_{G}(\mu)$ for some $H,G\in \mathcal{H}^\XX$, then
$$
\mu(H\Delta G)=0.
$$
\end{enumerate}
\end{theorem}
Theorem \ref{t2} is proved in Section \ref{s3}.

\begin{remark}\label{r3} Condition (ii) can be reformulated in the following, more geometrical, way. Without loss of generality, we can assume $\mu$ to be centered. Denote by $\mathcal{L}(\mu)$ the linear span of $D_0(\mu)$. By  statements (e), (f) of Theorem \ref{t1} and by the definition  of $Z(\mu)$, one has
\be\label{eL} \mathcal{L}(\mu)=\bigcup_{r\in
(0,\infty)}rD_{1/2}(\mu)=\bigcup_{r\in (0,\infty)}rZ(\mu). \ee Then
(ii) is equivalent to the following:
$$
\mathcal{L}(\mu)=\XX.\leqno(ii')
$$
By the identity (\ref{eL}), the linear space $\mathcal{L}(\mu)$ is
naturally endowed by the norm $\|\cdot\|_{\Lf}$  equal to  the
Minkowskii functional of the zonoid $Z(\mu)$; note that $Z(\mu)$ is
symmetric and convex because $\mu$ is centered. We call the normed
space $\mathcal{L}(\mu)$ with the norm $\|\cdot\|_{\Lf}$ \emph{the
zonoid span} of the measure $\mu$, and remark that condition (ii),
in fact, means that the initial space $\XX$ coincides with the
zonoid span of $\mu$.
\end{remark}

\begin{remark}
In the case  $\XX=\Re^d$ condition (ii$'$) (and therefore (ii))
holds trivially:  if (ii$'$) fails then the measure is concentrated
on the hyper-plane $\mathcal{L}(\mu)$, which contradicts to (i).

However, in the  infinite-dimensional setting
 condition (ii) is much more demanding.  For instance, this condition fails for any $\sigma$-additive probability measure $\mu$ on, say, an infinite-dimensional Hilbert space, such that the respective family (\ref{fam}) is uniformly integrable. Indeed, in that case  $D_0(\mu)$ is a  $\sigma$-compact set, and therefore by the Baire category theorem has an empty interior. In other words, in that case the geometry of the zonoid span $\mathcal{L}(\mu)$ necessarily differs from the geometry of $\XX$.
\end{remark}

\section{Proofs}\label{s3} In this section we prove Theorem \ref{t1}, Proposition \ref{p1},  and Theorem \ref{t2}. Without loss of generality we assume $\mu$ to be centered.

\subsection{Properties of zonoid and lift zonoid} Here  we prove statements (a) -- (c) of Theorem \ref{t1},  Proposition \ref{p1}, and Proposition \ref{p11}.

\emph{Statements (a), (b).} Clearly, it is enough to prove (a) and
(b) for the zonoid $Z(\mu)$, only. Convexity of $Z(\mu)$ follows by
the definition: if $x_1=E[\eta_1X]$ and $x_2=E[\eta_2X]$ are some
points from $Z(\mu)$, then for any $c\in [0,1]$ the point
$cx_1+(1-c)x_2=E[(c_1\eta_1+(1-c)\eta_2)X]$ belongs to $Z(\mu)$, as
well. Symmetricity also follows by the definition: for any $x=E[\eta
X]\in Z(\mu)$ the symmetric point $(-x)$ can be represented as
$(-x)=E[(1-\eta) X]$, and therefore belongs to $Z(\mu)$.

 Because for any $x=E[\eta X]\in Z(\mu)$ one has by (\ref{type1})
$$
|\la x, x^*\ra|\leq E|\eta||\la  X, x^*\ra|\leq E|\la  X,
x^*\ra|\leq C_\mu\|x\|^*_{\XX^*}, \quad x^*\in \XX^*,
$$
the set $Z(\mu)$ lies in the ball $\{x:\|x\|_\XX\leq C_\mu\}$ and
hence is bounded.

Let $x_n=E[\eta_n X], n\geq 1$ be a sequence in $Z(\mu)$. Consider
the variables $\eta_n, n\geq 1$ as elements of the unit ball in the
space $L_\infty(\Omega, \sigma(X), P)=(L_1(\Omega, \sigma(X),
P))^*$. Then by the Banach-Alaoglu theorem (see \cite{Rud73},
Chapter 3) there exists a $*$-weakly convergent subsequence
$\{\eta_{n_k}\}$, denote by $\hat \eta$ respective limit. Then
$$
\la x_{n_k}, x^*\ra=E\eta_{n_k}\la X, x^*\ra\to E\hat \eta \la X,
x^*\ra, \quad x^*\in \XX^*.
$$
Therefore,  if $x_n\to x$ in $\XX$, then $x=E[\hat \eta X]\in
Z(\mu)$; that is, $Z(\mu)$ is closed.\qed

\emph{Proposition 1.}  In fact, we have just proved that any
sequence $\{x_n\}\subset Z(\mu)$ has a subsequence $\{x_{n_k}\}$,
weakly convergent to $\hat x:=E[\hat \eta X]$. Let us prove that
$x_{n_k}\to \hat x$ (strongly) in $\XX$; this would give the
required compactness property.

Assuming that $x_{n_k}\not \to \hat x$ in $\XX$, we have a sequence
$\{x^*_k\}$ of elements of the unit ball in $\XX^*$ such that
\be\label{assume} \la x_{n_k}-\hat x, x^*_k\ra \not \to 0. \ee By
the  Banach-Alaoglu theorem, we can assume without loss of
generality that $x^*_k\to x^*$ $*$-weakly. Then
$$
\la X, x_k^*\ra \to \la X, x^*\ra,
$$
with probability 1. By the uniform integrability of (\ref{fam}),
this convergence also holds true in $L_1(\Omega, \sigma(X), P)$.
Then
$$
\la x_{n_k}, x_k^*\ra=E\eta_{n_k}\la X, x^*\ra+E\eta_{n_k}\Big(\la
X, x^*_k\ra-\la X, x^*\ra\Big)\to E\hat \eta \la X, x^*\ra=\la \hat
x, x^*\ra
$$
because $\eta_{n_k}\to\hat \eta$ $*$-weakly in $L_\infty(\Omega,
\sigma(X), P)$. The latter relation contradicts to (\ref{assume}).
This proves compactness of $Z(\mu)$. For $\hat Z(\mu)$ and
$D_\alpha(\mu), \alpha\in (0,1]$ the argument is completely
analogous.\qed

\emph{Statement (c).} The proof here is completely analogous to the
one of Theorem 3.5, \cite{KM_Be} in the finite-dimensional setting,
and hence we just sketch it. For $x^*\in X^*$, denote by $\mu_{x^*}$
respective one-dimensional projection of $\mu$; that is, the law of
the random variable $\la X, x^*\ra$. The calculation from the proof
of Theorem 3.4 in \cite{KM_Be} shows that for arbitrary cylinder
measures  $\mu, \nu$ the inclusion $\hat Z(\mu)\subset \hat Z(\nu)$
yields the inclusion $\hat Z(\mu_{x^*})\subset \hat Z(\nu_{x^*}),
x^*\in \XX^*$. Therefore, it follows from the identity $\hat Z(\mu)=
\hat Z(\nu)$ that for every $x^*\in \XX^*$ lift zonoids $\hat
Z(\mu_{x^*}), \hat Z(\nu_{x^*})$ coincide and hence (see Remark 3.2
in \cite{KM_Be}) $\mu_{x^*}=\nu_{x^*}$. Because a cylinder measure
is identified uniquely by its one-dimensional projections, we have
then $\mu=\nu$. \qed

\emph{Proposition 2.} We have just proved that $Z(\mu)$ is a
(non-empty) closed convex set. Then by by the second separability
theorem (e.g. \cite{Schaefer}, Chapter II, Theorem 9.2)
$$
Z(\mu)=\bigcap_{H\in \mathcal{H}^\XX: Z(\mu)\subset H} H.
$$
Take $x^*\in \SS_{\XX^*}$ and consider half-spaces of the form
$H=\{x:\la x, x^*\ra\leq a\}$. For every $\sigma (X)$-measurable
$\eta:\Omega\to [0,1]$, one has \be\label{incl} E\eta\la X,
x^*\ra\leq E\la X, x^*\ra_+=E h([0,X], x^*), \ee which means that
$Z(\mu)\subset H$ for any half-space of the above form with $a\geq E
h([0,X], x^*)$. On the other hand, for $\eta=\1_{\la X, x^*\ra\geq
0}$ inequality (\ref{incl}) turns into the equality, which means
that $Z(\mu)\not \subset H$ for any half-space of the above form
with $a< E h([0,X], x^*)$. This proves statement 1. Statement 2
follows from statement 1. \qed

\subsection{Properties of zonoid trimmed regions} Here  we prove statements (d) -- (h) of Theorem \ref{t1}. Statement (d) follows directly from statement (a), because $D_\alpha(\mu)$ is the projection on $\XX$ of the section of $\hat Z(\mu)$ by the affine hyper-plane $\{(t,x):t=\alpha\}\subset \Re\times \XX$. Statements (e) and (f) follow directly from the definition.

\emph{Statement (g).} By the statement (f), to prove the first part
of the statement it is sufficient to show that $$
D_{\alpha_n}(\mu)\to D_\alpha(\mu), \quad \alpha_n\uparrow \alpha<1
$$
in the Hausdorff distance. By the statement (e) we have
$D_\alpha(\mu)\subset D_{\alpha_n}(\mu)$, and the only thing we need
to prove is that for arbitrary sequence $x_n\in D_{\alpha_n}(\mu),
n\geq 1$ \be\label{dist} \mathrm{dist}\,(x_n, D_\alpha(\mu))\to 0.
\ee Let
$$
x_n={1\over \alpha_n}E[\eta_n X], \quad n\geq 1,
$$
where every $\eta_n:\Omega\to [0,1]$ is a $\sigma(X)$-measurable
random variable such that $E\eta_n=\alpha_n$. We put
$$
\zeta_{n}=\eta_{n}+ \frac{\alpha-\alpha_{n}}{1-
\alpha_{n}}(1-\eta_{n}),
$$
then $\zeta_n$ is a $\sigma(X)$-measurable random variable taking
values in $[0,1]$. Clearly, $E\zeta_n=\alpha$,  and  therefore
$x_n':=(1/\alpha)E[\zeta_nX]$ belongs to $D_\alpha(\mu)$. Because
 $\mu$ is centered, we have
$$x_n'=x_n-
\frac{\alpha-\alpha_{n}}{1- \alpha_{n}}x_n.
$$
Not that every $x_n$ belongs to the set $D_{\alpha_1}(\mu)$, which
is bounded by statement (d). Because $\alpha_n\to \alpha$ and
$\alpha<1$, we have $\|x_n-x_n'\|\to 0$, which proves (\ref{dist}).

The set $D_1(\mu)$ consists of one point $E[X]=0$. Hence, to prove
the second part of the statement,  we need to show that if
$x=(1/\alpha_n)E[\eta_n X]$ for some sequence of
$\sigma(X)$-measurable random variables with $\alpha_n:=E\eta_n\to
1$, then $x=0$. Since $\mu$ is centered,
$(-x)=(1/\alpha_n)E[(1-\eta_n)X]$. Clearly, $(1-\eta_n)\to 0$ in
probability. Then by the Lebesque dominated convergence theorem we
have
$$
E(1-\eta_n)\la X, x^*\ra\to 0,  \quad x^*\in \XX^*.
$$
This gives $\la(-x), x^*\ra=0, x^*\in \XX^*$ and therefore
$x=0$.\qed

\emph{Statement (h).}  First, let us show that for every $\alpha\in
(0,1]$ the set $D_\alpha(\mu)$ is contained in every $H\in
\mathcal{H}^\XX$ with $\mu(H)=1$; this would yield
$D_0(\mu)\subset\mathbf{S}(\mu)$. If $H=\XX$, the inclusion is
trivial. Otherwise $H=\{x:\la x, x^*\ra\geq a\}$ is an affine
half-space, and we have $\la X, x^*\ra\geq a$ with probability 1.
Then for every random variable taking values in $[0,1]$ one has
$$
E[\eta\la X, x^*\ra]\geq a E \eta,
$$
which implies the required inclusion.

 On the other hand, the  closure of $D_0(\mu)$ is convex by statements (d), (e) and, of course, is closed. Then, by the second separability theorem, for any point $x$ which does not belong to this closure there exist $x^*\in \XX^*, a\in \Re$ such that
$$
a:=\inf_{y\in D_0(\mu)}\la y, x^*\ra>\la x, x^*\ra.
$$
It is easy to see that the affine half-space $H=\{x:\la x,
x^*\ra\geq a\}$ satisfies $\mu(H)=1$. Indeed, if this is incorrect
then $\XX\setminus H=\{x:\la x, x^*\ra< a\}$ has non-zero measure
$\mu$, and then the point
$$z:={1\over\mu(\XX\setminus H)}E[\1_{\XX\setminus H}(X) X]$$
 belongs to $D_0(\mu)$. Since by the construction $\la z, x^*\ra<a$, this would contradict to the definition of $a$. Therefore, for arbitrary $x\not\in H$ there exists some  affine half-space $H$ with $\mu(H)=1$ such that $x\not \in H$. This means that  $x\not\in \mathbf{S}(\mu)$, and consequently the closure of $D_0(\mu)$ coincides with whole $\mathbf{S}(\mu)$.\qed

\subsection{Barycentric representation}\label{s33}    In this section, we use the properties (d) -- (h) of the family of zonoid $\alpha$-trimmed regions  to prove Theorem \ref{t2}.

First, we note that every $D_\alpha(\mu), \alpha\in (0,1)$ has a
non-empty interior. Indeed, it follows from condition (ii) and the
Baire category theorem that \emph{some} $D_\alpha(\mu), \alpha\in
(0,1)$ has a non-empty interior. Then the required statement follows
from the property (f).

Next, let us proceed with the construction of the barycentric
representation of the points  $x\in D_0(\mu)$ in the form
(\ref{bary}).
  Clearly, for  $x=0$ there exists  such representation with $H=\XX$, and this representation is unique in the sense explained above. Consider arbitrary  $x\in D_0(\mu)$, $x\not=0$. There exists  unique $\alpha=\alpha(x)\in (0,1)$ such that
$x\in D_\alpha(\mu)$ and $x\not \in D_{\beta}(\mu)$ for any
$\beta>\alpha$; this follows from the property (e) in Theorem
\ref{t1}.  By the continuity property (g), this yields $x\in
\partial D_\alpha(\mu)$. Because $D_\alpha(\mu)$ has non-empty
interior, we can apply the first separability theorem (e.g.
\cite{Schaefer}, Chapter II, Theorem 9.1 and subsequent corollary),
and get that there exists a \emph{support hyper-plane} at the point
$x$  to $D_\alpha(\mu)$. That is, there exists  $x^*\in
\XX^*\setminus \{0\}$ such that \be\label{compare2} \la x,
x^*\ra\geq\la y, x^*\ra, \quad y\in D_\alpha(\mu). \ee By the
condition (\ref{lines}), there exists $a$ such that
$$
\mu (\{x': \la x',x^*\ra \geq a \})=\alpha.
$$
 We put $f=\1_A(X), A=\{x':\la x',x^*\ra \geq a\}.$ Then $f\in F_\alpha(\mu)$, and
 for every $g\in F_\alpha(\mu)$
\be\label{compare} \ba E(f-g)\la X, x^*\ra&=E(f-g)\Big(\la X,
x^*\ra-a\Big)\\&=E\Big(\la X,
x^*\ra-a\Big)\Big((1-g)\1_A(X)\Big)+E\Big(\la X,
x^*\ra-a\Big)\Big(-g\1_{\XX\setminus A}(X)\Big)\geq 0, \ea \ee
because both expectations in the right hand side of (\ref{compare})
are non-negative. This means that
$$
{1\over \alpha}E[\1_{A}(X)X]=B_A(\mu)
$$
is the \emph{unique} point in the set $D_\alpha(\mu)$ which provides
the extremum for the functional $y\mapsto \la y, x^*\ra$ at this
set. Combined with (\ref{compare2}), this gives the required
representation (\ref{bary}) with $H=A$. The uniqueness of the
representation follows from the  strict inequality in
(\ref{compare}), valid as soon as $f$ and $g$ do not coincide a.s.

Note that, obviously, $\mathbf{B}(\mu)\subset D_0(\mu)$. Because we
have already proved the (unique) representation of any point of
$D_0(\mu)$ as a barycenter of some set $H\in \mathcal{H}^\XX$, this
completes the proof of   statement 2.

Recall that $D_0(\mu)$ is convex and has a non-empty interior. Then
(see \cite{Schaefer}, Chapter II, \S 1.3) the interior of $D_0(\mu)$
coincides with the interior of its closure. Combined with the
property (h) from Theorem \ref{t1}, this gives that every interior
point of $\mathbf{S}(\mu)$ belongs to (the interior of) $D_0(\mu)$.
On the other hand, since the interior of $\mathbf{S}(\mu)$ is
non-empty, it follows from the first separability theorem that for
every $x\in \mathbf{S}(\mu)$, which is not an interior point of
$\mathbf{S}(\mu)$, there exist $x^*\in \XX\setminus \{0\}$ and $a$
such that $\la x, x^*\ra=a$ and $\mu(x':\la x', x^*\ra\leq a)=1$.
Then, using (\ref{lines}), it is easy to show that for every every
point $y\in D_0(\mu)$ one has $\la y, x^*\ra<a$, and consequently
$x\not \in D_0(\mu)$. This means that $D_0(\mu)$ coincides with the
interior of $\mathbf{S}(\mu)$, and completes the proof of
statement 1.  \qed

\section{Barycentric coordinates}\label{s4}
We have mentioned briefly in the Introduction that the  barycentric
representation, established in Theorem \ref{t2}, in some particularly
interesting cases, makes it possible   to  introduce the
\emph{barycentric coordinates} on  the initial space. Here we
discuss this topic in details.

Assume that $\mu$ is a cylinder probability measure of type 1 on
$\XX$ which is centered and satisfies (\ref{lines}). Assume also
that \be\label{support} \mu(H\Delta G)>0\quad \hbox{for every $H,
G\in \mathcal{H}^\XX, \quad H\not=G$}. \ee Note that the above
assumptions yield that every non-empty set $H\in \mathcal{H}^\XX$
has positive measure $\mu$, and therefore  the set $\mathbf{S}(\mu)$
coincides with whole $\XX$.

 For $x\in \XX$ define
$$
\alpha(x)=\sup\{\alpha\in(0,1]:x\in D_\alpha(\mu)\},
$$
the \emph{zonoid data depth} of the point $x$ w.r.t. $\mu$, see
Definition 7.1 in \cite{KM_AnnSta}. For one exceptional point
$x=0$, we have $\alpha(x)=1$. For all others, we have $\alpha(x)\in
(0,1)$, and
 according to (\ref{support}) and the proof from Section \ref{s33} there exist unique $x^*\in \SS_{\XX^*}$  and $a\in \Re$ such that $x$ equals the barycenter of the half-space $H=\{x':\la x',x^*\ra \geq a\}$. Clearly, the pair $(a, x^*)$ determines $x$ uniquely, hence we can identify any point $x\in \XX\setminus\{0\}$ by a pair
\be\label{rep_a}
 (a, x^*) \in \Re\times \SS_{\XX^*}.
\ee
 On the other hand, for $x\not=0$ consider respective $x^*\in \SS_{\XX^*}$ and the value of the support function of  $D_{\alpha(x)}(\mu)$ on this $x^*$:
 $$
 h=h( D_{\alpha(x)}(\mu), x^*).
 $$
 It can be seen easily that $x$ is the unique point of  $D_{\alpha(x)}(\mu)$ on the affine hyper-plane $\{x':\la x', x^*\ra=h\}$. Hence, there is another possibility to  identify a point $x\in \XX\setminus\{0\}$ by a pair
 \be\label{rep_h}
 (h, x^*) \in (0,\infty)\times \SS_{\XX^*}.
 \ee
 Finally, given a set $D_\alpha(\mu)$ and $x^*\in \SS_{\XX^*}$, one can uniquely define $h>0$ such that the affine hyper-plane $\{x':\la x', x^*\ra=h\}$ is a support hyper-plane for $D_\alpha(\mu)$. Hence, we can also identify  a point $x\in \XX\setminus\{0\}$ by a pair
\be\label{rep_alpha}
 (\alpha, x^*) \in (0,1)\times \SS_{\XX^*},
 \ee
 where the first coordinate coincides with the zonoid data depth $\alpha(x)$.

Let us summarise. We have represented  $\XX$ as a disjunctive union
$$
\bigsqcup_{\alpha\in (0,1]}Z_\alpha(\mu),
$$
where $Z_\alpha(\mu)$ denotes the set of points with  the  zonoid
data depth w.r.t. $\mu$ equal $\alpha$. The set $Z_1(\mu)$ consists
of one point $0$. For $\alpha\in (0,1)$, any  point $x\in
Z_\alpha(\mu)$ is identified by respective $x^*\in \SS_{\XX^*}$,
which can be naturally understood as the unique  ``outside tangent
direction'' for $Z_\alpha(\mu)$ at the point $x$. This leads to the
representation (\ref{rep_alpha}). Without loss of the uniqueness of
representation property, the scalar $\alpha$ can be replaced either
by $h=\la x, x^*\ra$ (which would lead to (\ref{rep_h})), or by $a$
such that $x$ is the barycenter of the affine half-space $H=\{x':\la
x',x^*\ra \geq a\}$ (this  would lead to (\ref{rep_a})).

\section{Barycentric representation on a space with a $\sigma$-additive measure}\label{s5}

Let $\XX$ be a separable (not necessarily reflexive) Banach space
and $\mu$ be a probability measure on the Borel $\sigma$-algebra in
$\XX$. Assume that $\mu$ has a weak order of some order $p>1$; that
is, (\ref{wo}) holds true. Then Theorem \ref{t2} is not applicable
because its condition (ii) fails. Here we discuss the structure of
the zonoid span $\mathcal{L}(\mu)$, and give a version of the
barycentric representation valid in that case. Without loss of
generality we assume $\mu$ to be centered.

To study the properties of the zonoid span $\mathcal{L}(\mu)$, it is
convenient to use the objects related to the notion of a
\emph{measurable linear functional}; e.g. \cite{Dorogovtsev2000}.
Let us recall briefly respective constructions.

 For $r\in [1,p]$, denote the canonical embedding of $\XX^*$ into $L_r(\XX, \mu)$ by $T_{r, \mu}$. Define the space $\HH_{r, \mu}(\XX)$ as the closure of the range of $T_{r, \mu}$; any element  in $\HH_{r, \mu}(\XX)$ is called an \emph{($L_r$-integrable) measurable linear functional}.

For any $g\in L_{\infty}(\XX, \mu)$ the integral
$$
\int_\XX g(x)x\mu(dx)\in \XX
$$
is well defined, see \cite{VahTarCho}, Chapter II, \S3. It can be
verified easily that the adjoint operator
$T_{1,\mu}^*:L_{\infty}(\XX, \mu)\to \XX^{**}$ has the form
$$
T^*_{1,\mu}g=\int_\XX g(x)x\mu(dx), \quad g\in L_{\infty}(\XX, \mu),
$$
and in particular takes its values in $\XX\subset \XX^{**}$.

\begin{prop}\label{p2}

\begin{enumerate}\item The set $\mathcal{L}(\mu)$ equals the range of $T^*_{1,\mu}$.
\item The norm $\|\cdot\|_{\Lf}$ is equivalent to the norm $\|\cdot\|_{\Rf}$ defined by
\be\label{R} \|x\|_{\Rf}=\inf_{g\in L_\infty(\XX, \mu):
T^*_{1,\mu}g=x}\|g\|_{L_\infty(\XX, \mu)}, \quad x\in
\mathcal{L}(\mu)=\Rf(T^*_{1,\mu}). \ee
\end{enumerate}
\end{prop}

\emph{The proof} follows directly from the definitions of the
respective objects, and hence we omit the details. \qed

The statement of Proposition \ref{p2} can be made even more precise.
Consider the dual space $[\HH_1(\XX, \mu)]^*$; that is, the
factor-space in $L_{\infty}(\XX, \mu)$ w.r.t. the equivalence
\be\label{equiv} g_1\sim g_2 \Leftrightarrow
\int_\XX(g_1(x)-g_2(x))\la x, x^*\ra\mu(dx), \quad x^*\in \XX^*, \ee
see \cite{Rud73}, Chapter 4. One can see easily that for any $g\in
L_\infty(\XX, \mu)$
$$
T^*_{1,\mu}g=0 \Leftrightarrow g\sim 0.
$$
This means that $T^*_{1,\mu}$ can be naturally understood as a
linear operator $[\HH_1(\XX, \mu)]^*\to \mathcal{L}(\mu)$. In
addition, the norm in the space $[\HH_1(\XX, \mu)]^*$ is defined by
the relation
$$
\|g\|_{[\HH_1(\XX, \mu)]^*}=\inf_{\tilde g\sim g}\|\tilde
g\|_{L_\infty(\XX, \mu)}.
$$
Comparing this with the expression in the right hand of (\ref{R}),
we obtain the following.

\begin{corollary} The zonoid span $\mathcal{L}(\mu)$ is a Banach space isomorphic to the dual space of the space of $\HH_1(\XX, \mu)$ of $L_1$-integrable measurable linear functionals. Respective isomorphism is provided by the operator $T^*_{1,\mu}$.
\end{corollary}

Now, we are ready to give a version of the barycentric
representation valid in the case under the consideration. Denote by
$\mathcal{H}^\XX_{1,\mu}$ the family of the sets of the form
$$
\{x: h(x)\geq a\}, \quad h\in \HH_{1, \mu}(\XX), \quad a\in \Re.
$$

\begin{theorem}\label{t3}  Assume that

\begin{itemize}\item[(a)] for every $h\in \HH_{1, \mu}(\XX)\setminus\{0\}, a\in \Re$
$$
\mu(\{x: h(x)= a\})=0;
$$
\item[(b)] for some $r>1$, $\HH_{1, \mu}(\XX)=\HH_{r, \mu}(\XX).$
\end{itemize}

Then the following holds true.

\begin{enumerate}\item The set $\mathbf{B}(\mu)=\{B_A(\mu), \mu(A)>0\}$ of all barycenters w.r.t. $\mu$ equals the interior of $\mathbf{S}(\mu)\cap \Lf(\mu)$ w.r.t. to the topology induced by the norm $\|\cdot\|_\Lf$.
\item Any point $x\in \mathbf{B}(\mu)$ has a representation in the form
\be\label{bary}x=B_H(\mu), \quad H\in \mathcal{H}^\XX_{1,\mu}, \quad
\mu(H)>0.\ee This representation is unique in the following sense:
if  $B_{H}(\mu)=B_{G}(\mu)$ for some $H,G\in
\mathcal{H}^\XX_{1,\mu}$, then
$$
\mu(H\Delta G)=0.
$$
\end{enumerate}
\end{theorem}

\emph{Proof.} The identity $\HH_{1, \mu}(\XX)=\HH_{r, \mu}(\XX)$
yields that the Banach spaces  $[\HH_{1, \mu}(\XX)]^*$ and $[\HH_{r,
\mu}(\XX)]^*$ are isometric. Denote $r'=r/(r-1)$, then $[\HH_{r,
\mu}(\XX)]^*$ is a factor-space in $L_{r'}(\XX, \mu)$ w.r.t. to the
equivalence (\ref{equiv}). Using the fact that $L_{r'}(\XX, \mu)$ is
reflexive, it is easy to show that every space $[\HH_{r,
\mu}(\XX)]^*$, $[\HH_{1, \mu}(\XX)]^*$, and $\mathcal{L}(\mu)$ is
reflexive, as well.

Denote $\tilde \XX=\mathcal{L}(\mu),$ then $\tilde \XX^*$ is
isometric to $\HH_{1, \mu}(\XX)$. Denote by $J:\HH_{1, \mu}(\XX)\to \tilde \XX^*$  respective isometry, and define on the probability space $(\Omega, \mathcal{F}, P)=(\XX, \mathcal{B}(\XX), \mu)$ a generalised random
element $\tilde X$ valued in $\tilde \XX$ by
$$
\la \tilde X, Jh\ra =h, \quad h\in \HH_{1, \mu}(\XX).
$$
By the construction,  the law $\tilde \mu$ of $\tilde X$ is a
cylinder measure of the type 1. In addition, $\sigma(\tilde X)$
coincides with $\mathcal{F}$, which implies
that the zonoid, the lift zonoid, and the zonoid $\alpha$-trimmed
regions for $\mu$ are equal to the images of the same objects for
$\tilde \mu$ under the natural embedding $\mathcal{L}(\mu)\to \XX$.
Therefore the required statements follows from   Theorem \ref{t2}
applied to $\tilde \mu$. \qed

\begin{remark} Condition (b) in the above proof is used to guarantee that $\mathcal{L}(\mu)$ is a \emph{reflexive} Banach space. This condition
holds true in particularly interesting cases, i.e. for a Gaussian
measure $\mu$ or, more generally, for a convex measure $\mu$; see
\cite{Borell74}. However, this is a non-trivial structural
assumption on $\mu$, and one can be interested in other conditions
which would provide $\mathcal{L}(\mu)$ to be reflexive. This is a
subject for a further investigation.
\end{remark}

\section{Example}

We illustrate the above results by an example of a Gaussian measure
$\mu$, where the calculations are most simple. The  results below
are not genuinely new; e.g. statement 1 in Proposition \ref{p4} is
already known,  see Theorem 4.3 in \cite{Borell09}.

 Let $\XX$ be a separable Hilbert space, and $\mu$ be a centered cylinder Gaussian measure on $\XX$ with the identity covariance operator. By the definition, this means that for the respective generalized random element $X$
 $$
 (X, x)_\XX\sim \mathcal{N}(0, \|x\|^2_\XX), \quad x\in \XX;
 $$
here  we identify $\XX^*$ with $\XX$ and write $(\cdot, \cdot)_\XX$
instead of $\la\cdot, \cdot\ra$. Clearly, $\mu$ satisfies
(\ref{lines}) and (\ref{support}).

The law of $X$ is \emph{rotationally invariant}; that is, $UX$ has
the same distribution with $X$ for any unitary operator $U:\XX\to
\XX$. Then it is easy to see that every $D_\alpha(\mu)$ is a closed
ball in $\XX$; denote the radius of this ball by $r(\alpha)$.  We
can determine this radius easily, using the arguments from   Section
\ref{s4}. Namely, for every point $x$ on the boundary of
$D_\alpha(\mu)$ respective ``outside tangent direction'' equals
$\|x\|_\XX^{-1} x$. Then $x$ can be uniquely represented as the
barycenter of the half-space \be\label{half-space}
H_x=\left\{x':\left(x', \|x\|_\XX^{-1} x\right)_\XX\geq
a(x)\right\}, \ee where the scalar $a=a(x)$ is determined by the
relation $\mu(H_{x})=\alpha$. Because
$$
\mu(H_{x})=P\left(\left(X, \|x\|_\XX^{-1} x\right)_\XX\geq
a(x)\right),
$$
and $\left(X, \|x\|_\XX^{-1} x\right)$ has a standard Gaussian
distribution, we have
$$
a(x)=\Phi^{-1}(1-\alpha), \quad \Phi(u):=\int_{-\infty}^u\phi(v)\,
dv, \quad \phi(u):={1\over \sqrt{2\pi}} e^{-u^2/2}.
$$
Then
$$\ba
r(\alpha)&=\|x\|_\XX=\left(x, \|x\|_\XX^{-1} x\right)={1\over \alpha}E\left(X, \|x\|_\XX^{-1} x\right)_\XX\1_{(X, \|x\|_\XX^{-1} x)_\XX\geq a(x)}\\
&={1\over \alpha \sqrt{2\pi}}\int_{a(x)}^\infty u e^{-u^2/2}\,
du={1\over\alpha} I(\alpha), \ea
$$
where $I:=\phi\circ \Phi^{-1}$ is the so-called \emph{Gauss
isoperimetric function;} here we have used that
$\Phi^{-1}(1-\alpha)=-\Phi^{-1}(\alpha)$ and that $\phi$ is even.
Remark  that then  for $x\not=0$ respective zonoid data depth equals
$r^{-1}(\|x\|_\XX)$, and hence
$$
a(x)=-\Phi^{-1}(r^{-1}(\|x\|_\XX))=-G^{-1}(\|x\|_\XX),\quad\hbox{where}\quad
G(u):={\phi(u)\over \Phi(u)};
$$
note that the function $G:\Re\to (0, \infty)$ is monotonous.
Therefore,  we have proved the following.

\begin{prop}\label{p3} Let $\mu$ be a centered cylinder Gaussian measure on $\XX$ with the identity covariance operator.

Then the following holds.
\begin{enumerate}\item Every zonoid $\alpha$-trimmed region $D_\alpha(\mu), \alpha\in (0,1)$ equals the ball with the radius $r(\alpha)=\alpha^{-1}I(\alpha)$.

\item Every point $x\not=0$ has a unique representation as a barycenter w.r.t. $\mu$ of the half-space (\ref{half-space}) with $a(x)=-G^{-1}(\|x\|_\XX)$.
\end{enumerate}
\end{prop}

The case of $\mu$ being a centered $\sigma$-finite Gaussian measure
on a separable Banach space $\XX$, can be reduced to the above, like
it was done in Section \ref{s5}. Because convergence in probability
of a Gaussian sequence yields $L_p$-convergence for every $p\geq 1$,
we have that all spaces $\HH_{p, \mu}(\XX), p\geq 1$ now coincide.
Respective space is called \emph{the reproducing kernel space} for
the Gaussian measure $\mu$; we denote it by $\HH_\mu$. Because
$\HH_{\mu}=\HH_{2, \mu}(\XX)$ can be considered as the Hilbert
space,  $\mathcal{L}(\mu)$ equals the image of $\HH_\mu$ under the
natural embedding $J_\mu:\HH_\mu\to \XX$ (this image is called
usually the \emph{Cameron-Martin space} of $\mu$). On this space,
$\mu$ generates a centered Gaussian cylinder measure $\tilde \mu$
with the identity covariance operator. Applying Proposition \ref{p3} to
$\tilde \mu$, we arrive at following.

\begin{prop}\label{p4} Let $\mu$ be a a centered $\sigma$-finite Gaussian measure on a separable Banach space $\XX$.

Then the following holds.
\begin{enumerate}\item Every zonoid $\alpha$-trimmed region $D_\alpha(\mu), \alpha\in (0,1)$ equals the image   of the ball in $\HH_\mu$ with the radius $r(\alpha)=\alpha^{-1}I(\alpha)$ under the natural embedding of $\HH_\mu$ into $\XX$.

\item Every point $x\not=0$ in the Cameron-Martin space has a unique representation as a barycenter w.r.t. $\mu$ of the measurable half-space $$
H_x=\left\{x':\left(x', \|J^{-1}_\mu x\|_{\HH_\mu}^{-1} J^{-1}_\mu
x\right)_{\HH_\mu}\geq a(x)\right\},\quad  a(x)=-G^{-1}(\|J^{-1}_\mu
x\|_{\HH_\mu}),
$$
where $(\cdot, h)_{\HH_\mu}$ denotes the linear measurable
functional respective to $h\in \HH_\mu$.
\end{enumerate}
\end{prop}

\textbf{Acknowledgement.} The authors a grateful to Ilya Molchanov,
who have turned their attention to the lift zonoid concept.


\end{document}